\documentclass[a4paper, 10pt, conference,twocolumn]{cssconf} 

\IEEEoverridecommandlockouts    
\usepackage{psfrag}
\usepackage{amsmath,amssymb,amsfonts,bbm}
\usepackage{latexsym}
\usepackage{graphicx}

\usepackage{colortbl}
\usepackage{fancyhdr}
\usepackage{epstopdf}
\usepackage{float}
\usepackage{hyperref}
\usepackage{booktabs}
\usepackage{enumerate}
\usepackage[ruled]{algorithm2e}
\usepackage{balance}
\usepackage{mathrsfs}

\usepackage[usenames,dvipsnames,svgnames]{xcolor}

\newtheorem{theorem}{Theorem}
\newtheorem{definition}{Definition}
\newtheorem{proposition}{Proposition}

\newtheorem{remark}{Remark}

\newtheorem{standing}{Standing Assumption}

\IEEEoverridecommandlockouts

\newcommand{\R}{\mathbb{R}}

\newcommand{\mc}{\mathcal}

\newcommand{\col}{\mathrm{col}}

\newcommand{\dom}{\mathrm{dom}}
\newcommand{\proj}{\mathrm{proj}}

\newcommand{\bs}{\boldsymbol}

\graphicspath{{Figures/} }

\newcommand{\Rmnum}[1]{\expandafter\@slowromancap\romannumeral #1@}

\newcommand{\dst}{\displaystyle}

\def\be{\begin{equation}}
\def\ee{\end{equation}}
\def\ba{\begin{array}}
\def\ea{\end{array}}
\def\eqa{\begin{eqnarray}}
\def\eqe{\end{eqnarray}}

\newtheorem{lem}{Lemma}

\usepackage[prependcaption,colorinlistoftodos]{todonotes}

\begin{document}

\title{Continuous-time integral dynamics for aggregative game equilibrium seeking}

\author{Claudio De Persis, Sergio Grammatico 
\thanks{C. De Persis is with the Faculty of Science and Engineering, University of Groningen, The Netherlands. S. Grammatico is with the Delft Center for Systems and Control (DCSC), TU Delft, The Netherlands. E-mail addresses: \texttt{c.de.persis@rug.nl}, \texttt{s.grammatico@tudelft.nl}. This work was partially supported by the Netherlands Organisation for Scientific Research (NWO) under research projects OMEGA (grant n. 613.001.702) and P2P-TALES (grant n. 647.003.003).
}
}

\thispagestyle{empty}
\pagestyle{empty}

\maketitle

\begin{abstract}
In this paper, we consider continuous-time semi-decentralized dynamics for the equilibrium computation in a class of aggregative games. Specifically, we propose a scheme where decentralized projected-gradient dynamics are driven by an integral control law. To prove global exponential convergence of the proposed dynamics to an aggregative equilibrium, we adopt a quadratic Lyapunov function argument. We derive a sufficient condition for global convergence that we position within the recent literature on aggregative games, and in particular we show that it improves on established results.
\end{abstract}

{\keywords
Noncooperative game theory, Multi-agent systems, Decentralized control, Projected dynamical systems.
}

\medskip
\section{Introduction}

Aggregative game theory \cite{kukushkin:04} represents a mathematical framework to model  inter-dependent optimization problems for a set of noncooperative agents, or players, whenever the decision of each agent is affected by some aggregate effect of all the agents. Remarkably, this feature emerges in several application areas, such as demand side management in the smart grid \cite{Saad2012}, e.g. for electric vehicles \cite{parise:colombino:grammatico:lygeros:14, ma:zou:ran:shi:hiskens:16, grammatico:16cdc-pev} and thermostatically controlled loads \cite{grammatico:gentile:parise:lygeros:15, Li2016}, demand response in competitive markets \cite{li:chen:dahleh:15} and network congestion control \cite{barrera:garcia:15, grammatico:17}.

Existence and uniqueness of equilibria in (aggregative) games has been widely studied in the literature \cite{facchinei2007finite, facchinei2010generalized}, \cite[\S12]{palomar2010convex}, \cite[Part II]{cominetti:facchinei:lasserre}, also in automatic control and control systems\cite{pavel:07, yin:shanbhag:mehta:11, kulkarni:shanbhag:12}. To compute a game equilibrium, severable algorithms are available, both distributed \cite{salehisadaghiani:pavel:16,koshal:nedic:shanbhag:16,ye:hu:17, ye:hu:17tac,liang:yi:hong:17} and semi-decentralized  \cite{grammatico:parise:colombino:lygeros:16,grammatico:16cdc-convex, grammatico:17, belgioioso:grammatico:17cdc,belgioioso:grammatico:18ecc}. 

The majority of the available algorithms consider discrete-time dynamics that, under appropriate technical assumptions and sufficient conditions on the problem data, converge to an equilibrium of the game, e.g. Nash equilibrium. An elegant approach is in fact to characterize the desired equilibrium solutions as the zeros of a monotone operator, e.g.\ the concatenation of interdependent Karush--Kuhn--Tucker operators, and in turn formulate an equivalent fixed-point problem, which is solved via appropriate discrete-time nonlinear dynamics with guaranteed global asymptotic convergence.

On the other hand, only a few contributions have addressed the equilibrium computation problem via continuous-time dynamics, e.g.\ \cite{ye:hu:17tac,liang:yi:hong:17}. One reason is that in the continuous-time case where optimal decision problems of the agents are subject to constraints, projected dynamics, hence dynamics with \textit{discontinuous} right-hand side, shall be considered. Specifically, in \cite{ye:hu:17tac}, the authors propose distributed dynamics based on a consensus protocol to compute a Nash equilibrium, under technical assumptions that include the strong monotonicity of the so-called pseudo-gradient game mapping, but in the absence of constraints. In \cite{liang:yi:hong:17}, the authors address the generalized Nash equilibrium problem, that is, with both local and coupling constraints, via extended continuous-time dynamics that include three auxiliary vector variables for each agent. To prove convergence to an equilibrium, the authors postulate that the pseudo-gradient mapping is strictly monotone and that the parameter gains are chosen large enough, that is, directly proportional to the number of agents in the game.

Differently from the available literature, and in particular from \cite{ye:hu:17tac,liang:yi:hong:17}, the contribution of this paper is to provide a first, simple, integral control algorithm for the computation of an aggregative equilibrium via semi-decentralized dynamics. Since we consider games with constraints, we propose equilibrium seeking dynamics that are described as a projected dynamical system \cite{nagurney:zhang}. Therefore, the technical difficulty we address is to study the solutions to the derived projected dynamical system in the sense of Carath\`{e}odory solutions. With this aim, we exploit an invariance principle for projected dynamical systems. The main technical contribution is then to prove global exponential convergence of the proposed continuous-time, semi-decentralized, dynamics to an equilibrium of the considered aggregative game, under a mild design choice that involves a few problem parameters only. Interestingly, we discover that the derived sufficient condition improves on that in \cite[Th. 2]{grammatico:16cdc-convex}, especially when the number of agents playing the aggregative game is large.

The paper is organized as follows. Section \ref{sec:aggregative.games} introduces the mathematical setup and defines the aggregative equilibrium as solution concept. In Section \ref{sec:integral.dynamics}, we propose integral dynamics for the equilibrium computation and present the main technical results. In Section \ref{sec:discussion}, we compare the derived sufficient condition for global convergence within the recent literature of (monotone) aggregative games. 
Section \ref{sec:conclusion} summarizes the message of the paper and points at future research directions. The proofs are given in the Appendix to improve the reading flow.

\newpage
\subsection*{Basic notation}
$\R$ denotes the set of real numbers, while $\overline{\R} := \R \cup \{\infty\}$ the set of extended real numbers.
$I$ denotes the identity matrix. $\bs{0}$ ($\bs{1}$) denotes a matrix/vector with all elements equal to $0$ ($1$); to improve clarity, we may add the dimension of these matrices/vectors as subscript. 
$A \otimes B$ denotes the Kronecker product between matrices $A$ and $B$. $\left\| A \right\|$ denotes the maximum singular value of matrix $A$.
Given $N$ vectors $x^1, \ldots, x^N \in \R^n$, we define $\boldsymbol{x} := \col\left(x^1,\ldots,x^N\right) = \left[ {x^1}^\top, \ldots , {x^N}^\top \right]^\top$, and $\boldsymbol{x}^{-i} := \col\left(x^1,\ldots , x^{i-1}, x^{i+1},\ldots,x^N\right)$. Given $N$ sets $\mathcal{X}^1, \ldots, \mathcal{X}^N$, we define the Cartesian product by $\bs{\mathcal{X}} := \mathcal{X}^1 \times \ldots \times \mathcal{X}^N$.

\subsection*{Definitions}
Let the set $\mathcal{S} \subseteq \R^n$ be non-empty. The mapping $\iota_{\mc{S}}:\R^n \rightarrow \{ 0, \, \infty \}$ denotes the indicator function, i.e., $\iota_{\mc{S}}(x) = 0$ if $x \in \mathcal{S}$, $\infty$ otherwise. The set-valued mapping $N_{\mc{S}}: \R^n \rightrightarrows \R^n$ denotes the normal cone operator, i.e., 
$N_{\mc{S}}(x) = \varnothing$ if $x \notin \mc{S}$, $\left\{ v \in \R^n \mid \sup_{z \in \mc{S}} \, v^\top (z-x) \leq 0  \right\}$ otherwise.
The mapping $\proj_{\mathcal{S}}(\cdot) := \textrm{argmin}_{ y \in \mathcal{S}} \left\| y - \cdot\right\| : \R^n \rightarrow \mathcal{S}$ denotes the projection operator; $\Pi_{ \mathcal{S} }(x,v) := \lim_{ \epsilon \rightarrow 0^{+} } \frac{ \proj_{ \mathcal{S}}( x + \epsilon v ) - x }{ \epsilon } $ denotes the projection of the vector $v \in \R^n$ onto the tangent cone of $\mathcal{S}$ at $x \in \mathcal{S}$.

For a function $f: \R^n \rightarrow \overline{\R}$, $\dom(f) := \{x \in \R^n \mid f(x) < \infty\}$; $\partial f: \dom(f) \rightrightarrows {\R}^n$ denotes its subdifferential set-valued mapping, defined as $\partial f(x) := \{ v \in \R^n \mid f(z) \geq f(x) + v^\top (z-x)  \textup{ for all } z \in \textup{dom}(f) \}$; if $f$ is differentiable at $x$, then $\partial f(x) = \left\{ \nabla f(x)\right\}$. 

Given a closed convex set $\mathcal{C} \subseteq \R^n$ and a single-valued mapping $\varphi: \mathcal{C} \rightarrow \R^n$, the variational inequality problem VI$(\mathcal{C},\varphi)$, is the problem to find $x^* \in \mathcal{C}$ such that 
$(y-x^*)^\top  \, \varphi(x^*) \geq 0 \ \textup{  for all } y\in \mathcal{C}.$

\medskip
\section{Technical background: Aggregative games and variational aggregative equilibria}
\label{sec:aggregative.games}

An aggregative game, $G_{\textup{a}}\left( \mathcal{I}, ( J^i )_{i\in \mathcal{I}}, (\mathcal{X}^i)_{i \in \mathcal{I}} \right)$, consists of $N$ agents (or players) indexed by the set $\mathcal{I}:=\{1,2,\ldots, N\}$, where each agent $i$ can decide on a vector $x^i \in \mathcal{X}^i \subseteq \mathbb{R}^{n}$, with the aim to minimize its cost function $\bs{x} \mapsto J^i (x^i, {\rm avg}(\bs{x}))$, $J^i: \R^n \times \R^n \rightarrow \R \cup \{\infty\}$, where ${\rm avg}( \bs{x} ) := \frac{1}{N}\sum_{j\in \mathcal{I}} x^j$.

Let us postulate the following standing assumption that is valid throughout the paper.

\smallskip
\begin{standing} \label{asspt:strong.convexity.aggr}
\textit{Compactness and convexity}, from~\cite[Assumptions 1, 2]{grammatico:16cdc-convex}.
For each $i \in \mathcal{I}$, the set $\mathcal{X}^i \subsetneq \mathbb{R}^{n}$ is nonempty, compact and convex; the cost function $J^i$ is defined as
\be\label{cost.aggregative}
J^i(x^i, \sigma) = f^i(x^i) + \left( C \sigma \right)^\top x^i + \iota_{\mathcal{X}^i}\!\left( x^i \right), 
\ee
where the function $f^i: \R^n \to \mathbb{R}$ is twice continuously differentiable and $\ell$-strongly convex in $\mathcal{X}^i$, with constant $\ell>0$, and $C \in \R^{n \times n}$.
{\hfill $\square$}
\end{standing}
\smallskip

\smallskip
\begin{remark}
In a discrete-time setting, the functions $\{f_i\}_{i \in \mathcal{I}}$ can be assumed to be continuous, not necessarily differentiable \cite{grammatico:16cdc-convex}, \cite{grammatico:17}. Since the focus of this paper is equilibrium seeking in continuous-time, to avoid technical difficulties with non-smooth analysis and set-valued subdifferential mappings, we assume that the functions are twice continuously differentiable, as in \cite{ye:hu:17tac,liang:yi:hong:17}. We envision, however, that the differentiability assumption can be relaxed via the arguments in \cite{goebel:17}. 
{\hfill $\square$}
\end{remark}
\smallskip

In this paper, we focus on the design of continuous-time dynamics that converge to an aggregative equilibrium, which is a set of decision variables such that each is optimal given the average among all the decision variables.

\smallskip
\begin{definition}
\textit{Aggregative equilibrium}, from \cite[Def. 1]{grammatico:17}.
A collective vector ${\rm col}\left(\bar{x}^{1} , \ldots , \bar{x}^{N} \right) = \bs{\bar{x}}\in \bs{\mathcal{X}}$ is an aggregative equilibrium for the game $G_{\textup{a}}\left( \mathcal{I}, ( J^i )_{i\in \mathcal{I}}, (\mathcal{X}^i)_{i \in \mathcal{I}} \right)$ if, for all $i \in \mathcal{I}$,  
$$
\bar{x}^{i} \in \underset{y\in \mathcal{X}^i}{\textrm{argmin}} \, J^i\left(y,{\rm avg}(\bs{\bar{x}})\right).
$$
The set of aggregative equilibria is denoted by $\bar{\bs{\mathcal{X}}}$.
{\hfill $\square$}
\end{definition}
\smallskip

\smallskip
\begin{remark}
Alternatively, instead of the aggregative equilibrium, the notion of Nash equilibrium can be considered as target collective set, see \cite[Sec. V]{grammatico:16cdc-convex}. Moreover, we refer the interested reader to \cite[Th. 1]{grammatico:17} and to \cite{vaya:grammatico:andersson:lygeros:15,deori:margellos:prandini:17} for a comparison between aggregative and Nash equilibria.
{\hfill $\square$}
\end{remark}
\smallskip

It follows from \cite[Prop. 1]{grammatico:17} that an aggregative equilibrium does exist.
In view of Standing Assumption \ref{asspt:strong.convexity.aggr}, we provide an equivalent characterization of aggregative equilibrium via a variational inequality.

\smallskip
\begin{lem}
\label{lem:cnes.agg.eq}
\textit{Variational characterization}, from \cite[Prop. 1.4.2]{facchinei:pang}.
A collective vector $\bs{\bar{x}} = {\rm col}\left(\bar{x}^1, \ldots, \bar{x}^N\right) \in \bs{\mathcal{X}}$ is an aggregative equilibrium for the aggregative game $G_{\textup{a}}\left( \mathcal{I}, ( J^i )_{i \in \mathcal{I}}, (\mathcal{X}^i)_{i \in \mathcal{I}} \right)$ if and only if it satisfies the variational inequality
\be\label{VI}
\inf_{ z \in \mathcal{X}^i } \left(z-\bar{x}^i\right)^\top \!\left( \partial f^i(\bar{x}^i)+ C \, {\rm avg}(\bs{\bar{x}})\right) \ge 0, \ \ \forall i\in \mathcal{I}.
\ee
\hfill $\square$
\end{lem}
\smallskip

In view of the equilibrium characterization in Lemma \ref{lem:cnes.agg.eq}, uniqueness of the aggregative equilibrium can be studied via variational inequality arguments, see \cite[Th. 3] {scutari:facchinei:pang:palomar:14}.

\medskip
\section{Continuous-time integral dynamics for aggregative equilibrium seeking}
\label{sec:integral.dynamics}

With the aim to asymptotically reach an aggregative equilibrium, we propose the continuous-time dynamics
\be\label{aggregative.dynamics}
\ba{rcll}
\dot x^i &=& \Pi_{\mathcal{X}^i} \left( x^i \, , \, - \partial f_i(x^i)-  C\sigma \right), & \forall i\in \mathcal{I}, 
\smallskip \\ 
\dot \sigma &=& k \left( -\sigma + \frac{1}{N} \sum_{j \in \mathcal{I}} x^j \right)
\ea
\ee
which reads in projected vector form as
\be\label{aggregative.dynamics.vector}
\begin{bmatrix} 
\dot{\bs{x}} \\ \dot \sigma \end{bmatrix} = 
\Pi_{\mathcal{X} \times \R^n } \left( 
\left[\begin{matrix}
\bs{x} \\
\sigma \end{matrix}\right]
 \, , \, \begin{bmatrix}
 - \partial f( \bs{x} ) - \bs{1}_N \otimes ( C \sigma ) \\
 k \left(-\sigma +{\rm avg}( \bs{x} ) \right)
\end{bmatrix}  \right)
\ee
where $\partial f (\bs{x}) = \col\left(\partial f^1(x^1), \ldots, \partial f^N (x^N)\right)$ and $k > 0$ is a free design parameter.

We note that the structure of the computation and information exchange in \eqref{aggregative.dynamics} is semi-decentralized: the agents perform decentralized computations, namely, projected-gradient steps, and do not exchange information among each other, while a central control unit, which does not participate in the game, collects the average of the agents' decisions (aggregate information), ${\rm avg}(\bs{x}(t))$, and broadcasts a signal, $\sigma(t)$, to the agents playing the aggregative game. In turn, the dynamics of the broadcast signal are driven by the aggregate information, ${\rm avg}(\bs{x}(t))$. 

First, we observe that an equilibrium for the dynamics in \eqref{aggregative.dynamics} generates an aggregative equilibrium.

\smallskip
\begin{lem}\label{lem:agg.eq.and.dyn.eq}
If the pair $\left( \bs{\bar x}, \bar \sigma \right)$ is an equilibrium for the dynamics in \eqref{aggregative.dynamics}, then $\bs{\bar x}$ is an aggregative equilibrium for the game $G_{\textup{a}}\left( N, \{ J^i \}_{i}, \{ \mathcal{X}^i \}_{i} \right)$. Conversely, if $\bs{\bar x}$ is an aggregative equilibrium for the game $G_{\textup{a}}\left( \mathcal{I}, ( J^i )_{i \in \mathcal{I}}, ( \mathcal{X}^i )_{i\in \mathcal{I} } \right)$, then the pair $\left( \bs{\bar x}, \bar \sigma \right) =  \left( \bs{\bar x}, {\rm avg}(\bs{\bar x}) \right)$ is an equilibrium for the dynamics in \eqref{aggregative.dynamics}.
\hfill $\square$
\end{lem}
\smallskip

In the rest of the section, we analyze the convergence of the projected dynamics in \eqref{aggregative.dynamics.vector} to an equilibrium. In particular, we take two preliminary steps. First, we introduce a quadratic storage function, $V(\bs{x})$, and analyze its Lyapunov derivative, $\nabla V (\bs{x}) \, \dot{\bs{x}}$. Then, we consider a quadratic Lyapunov function candidate, $W\left( \bs{x}, \sigma \right)$, and establish a condition on the problem data $\ell$, $C$, $N$ and on the design parameter $k$ such that its Lyapunov derivative is negative definite. Finally, we establish our main asymptotic convergence result, whose proof - given in the appendix - is based on {Lyapunov stability theory} for projected dynamical systems with Carath\`{e}odory solutions.

\smallskip
\begin{lem} \label{cor.key}
\textit{Storage Lyapunov function}.
The function 
$$ \bs{x} \mapsto V( \bs{x} ) := \tfrac{1}{2} \| \bs{x} - \proj_{\bar{\bs{\mathcal{X}}}}( \bs{x} ) \|^2$$ 
is such that, for all  $x \in \bs{\mathcal{X}}$ and $\bs{u} \in \R^{n N}$, 
\begin{multline}
\dst {\partial V}(\bs{x})^\top \, \Pi_{\bs{\mathcal{X}}} \left( \bs{x} \, , \, - {\partial f}( \bs{x} ) + \bs{u} \right) \leq \\
-\left(\bs{x}-\proj_{\bar{\bs{\mathcal{X}}}}(\bs{x})\right)^\top \left({\partial f}(\bs{x})-{\partial f}(\proj_{\bar{\bs{\mathcal{X}}}}(\bs{x})) + \bs{u} - \bar{\bs{u}} \right),
\end{multline}
where $\bar{\bs{u}} := -\bs{1}_N \otimes {\rm avg} \left(\proj_{\bar{\bs{\mathcal{X}}}}(\bs{x})\right)$.
\hfill $\square$
\end{lem}
\smallskip

\smallskip
\begin{proposition} \label{prop:lyapunov.decrease}
\textit{Lyapunov decrease}.
Consider the Lyapunov function candidate $W: \bs{\mathcal{X}} \times \mathbb{R}^n \rightarrow \R$ defined as 
\[
W(\bs{x},\sigma)= \tfrac{1}{2} \| \bs{x} -\proj_{\bar{\bs{\mathcal{X}}}}(\bs{x})\|^2 + \tfrac{1}{2} \|\sigma - {\rm avg}( \proj_{\bar{\bs{\mathcal{X}}}}(\bs{x}) ) \|^2.
\]
If the inequality condition 
\be\label{stab.cond}
\min \{ \ell, k \} > \textstyle \tfrac{1}{2} \|C\|_\infty + \tfrac{1}{2} \frac{k}{N},
\ee 
holds, then we have that 
\begin{equation}
\label{strict.lyap}
\dot W(\bs{x},\sigma) := \nabla W(\bs{x},\sigma) 
\left[ 
\begin{smallmatrix}
\dot{\bs{x}} \\
\dot{\sigma}
\end{smallmatrix}
\right] \leq -\epsilon \, W(\bs{x},\sigma)
\end{equation}
for all $(\bs{x},\sigma)\in \bs{\mathcal{X}} \times \mathbb{R}^n$, for some $\epsilon > 0$. 
\hfill $\square$
\end{proposition}
\smallskip

\smallskip
\begin{theorem} \label{main.aggreg}
\textit{Global exponential convergence}.
Assume that the condition in \eqref{stab.cond} holds. For any initial condition $\left( \bs{x}_0, \sigma_0 \right) \in \bs{\mathcal{X}} \times \mathbb{R}^n$, there exists a unique Carath\`{e}odory solution to  \eqref{aggregative.dynamics.vector} starting from $\left( \bs{x}_0, \sigma_0 \right)$, which remains in $ \bs{\mathcal{X}} \times \mathbb{R}^n$ for all time, and exponentially converges to the set $\bar{\bs{\mathcal{X}}} \times {\rm avg}\left( \bar{\bs{\mathcal{X}}} \right)$.
\hfill $\square$
\end{theorem}

\medskip
\section{Discussion on the convergence condition} \label{sec:discussion}

For $k=1$, the $\sigma$-dynamics in \eqref{aggregative.dynamics.vector} can be thought as a control law, implemented by a central agent, which is the continuous-time counterpart of the Banach--Picard iteration studied in \cite{grammatico:16cdc-convex}, while for $k\in (0,1)$, the counterpart of the Krasnoselskij iteration.

We note however that in the continuous-time dynamics in \eqref{aggregative.dynamics.vector}, the parameter gain $k$ must not necessarily be smaller than $1$. 
For instance, if $k = \alpha \, \ell$, for some $\alpha > 0$, then the inequality condition in 
\eqref{stab.cond} 
reads as $2 \ell > \| C\|_\infty +\frac{\alpha \ell}{N}$, hence $\left( 2 - \tfrac{\alpha}{N}\right) \ell > \| C\|_\infty$, which is equivalent to $\ell > \frac{N}{2N-\alpha} \| C\|_\infty$ whenever $\alpha < 2N$. Let us compare the latter with the inequality condition established in  \cite[Th. 2]{grammatico:16cdc-convex}, i.e., $\ell \ge \| C\|$. Since $\| C\|_\infty \leq \| C\|$, the inequality condition in \eqref{stab.cond} is less strict than that in \cite[Th. 2]{grammatico:16cdc-convex} if $\frac{N}{2N-\alpha} < 1$, i.e., $\alpha < N$, e.g.\ for large number of agents $N$.

Next, we relate the condition in \eqref{stab.cond} with the monotonicity of the mapping that defines the variational inequality in \eqref{VI}, Lemma \ref{lem:cnes.agg.eq}, which is 
\begin{equation}
\label{eq:F}
\bs{x} \mapsto F(\bs{x}) := \partial f(\bs{x}) + \bs{1}_N \otimes \left( C {\rm avg}(\bs{x}) \right).
\end{equation}
In view of \cite[Cor. 1]{belgioioso:grammatico:17} and \cite[Lemma 3]{grammatico:parise:colombino:lygeros:16}, since the local cost functions $\{f^i\}_{i \in \mathcal{I}}$ are $\ell$-strongly convex, we have that the mapping $F$ is strictly monotone if $\ell + \tfrac{1}{2}\lambda_{\min}\left( C + C^\top\right) > 0$, e.g.\ if $C + C^\top \succcurlyeq 0$.

Finally, differently from the convergence condition provided in \cite[Sec. 4]{liang:yi:hong:17}, where the parameter gains must be chosen directly proportional to the number of agents, $N$, the inequality condition in \eqref{stab.cond} becomes less strict as $N$ grows. Clearly, the latter is an advantageous feature especially when the number of agents is very large.

\medskip
\section{Numerical simulation: Illustrative game-theoretic demand side management} \label{sec:applicability}

A recurrent application area for aggregative games is the demand side management, where the agents of the game are electricity prosumers, that is, consumers and possibly producers, whose disutility function depends on the average electricity demand, e.g.\ via a common electricity price.

For instance, in \cite[Sec. 5.2]{liang:yi:hong:17}, the individual cost function of each agent $i$ is chosen as $$ J^i\left( x^i, \sigma \right) = \tfrac{1}{2} \ell_i \left\| x^i - x_{*}^i \right\|^2 + \left( a \, \sigma + b \right)^\top x^i, $$ for some reference decision variable $x_{*}^i \in \R^n$, parameter $\ell_i > 0$, parameters $a \in \R$ and $b \in \R^n$. Since the dependence of $J^i$ on its second argument is affine, the setup satisfies our Standing Assumption \ref{asspt:strong.convexity.aggr}, precisely with $C = a \, I_n$.

\begin{figure}[!h]
\begin{center}
\includegraphics[width=1\columnwidth]{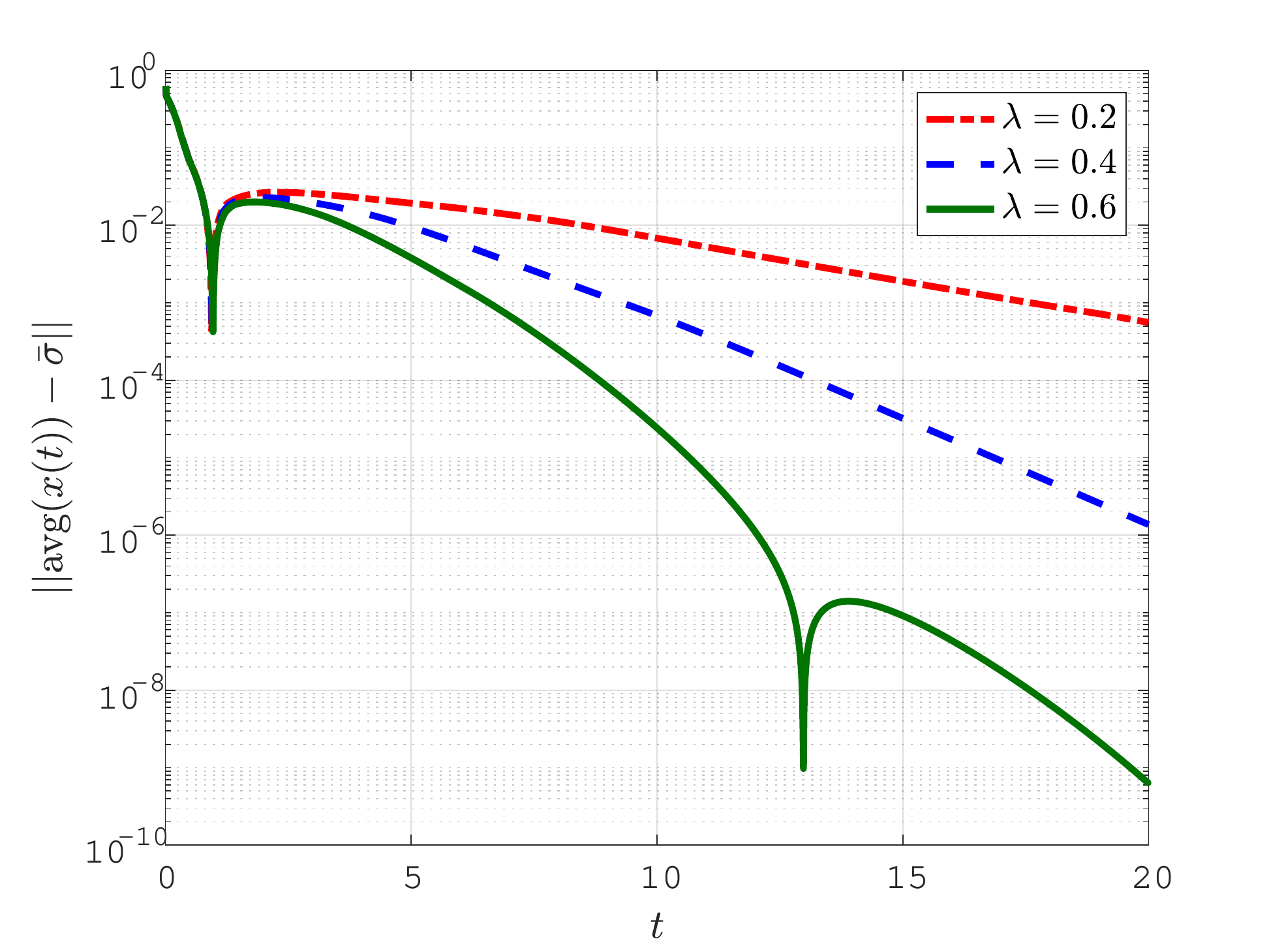}
\caption{Distance between the average strategy over time and the average strategy at the aggregative equilibrium.}
\label{fig:avg}
\end{center}
\end{figure}

\begin{figure}[!h]
\begin{center}
\includegraphics[width=1\columnwidth]{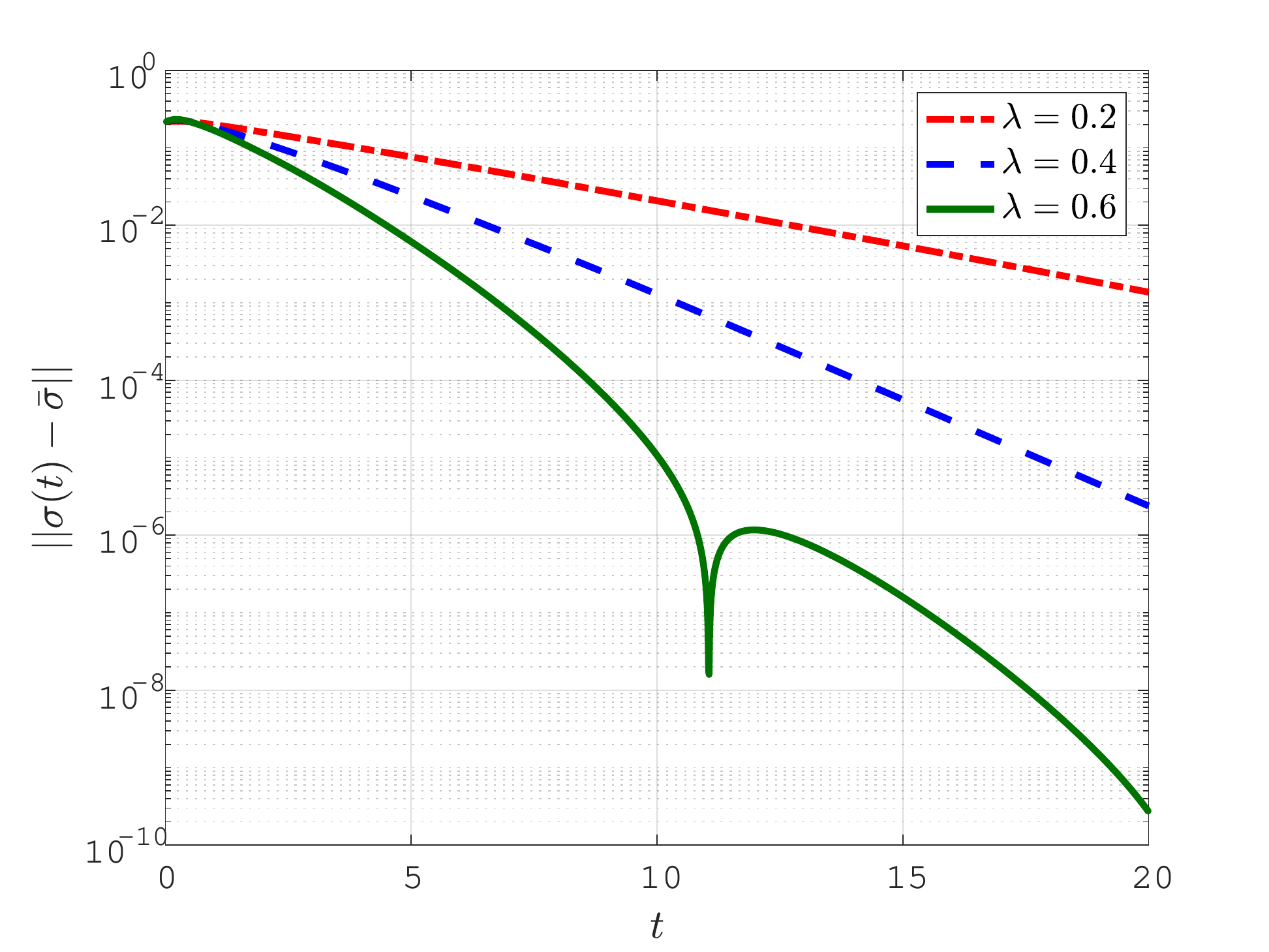}
\caption{Distance between the integral control signal over time and the average strategy at the aggregative equilibrium.}
\label{fig:sigma}
\end{center}
\end{figure}

Numerically, we take the following values: $N=100$, $a=1$, $b=0.5$, $\mathcal{X}^i = \left[ 0.25 \, , \, 0.75\right]$, $x_{*}^{i}$ sampled uniformly in $[0,1]$ and $\ell_i = 1.5$ for all $i$. 
For three choices of $k$, i.e., $k=0.2$, $k=0.4$, or $k=0.6$, we plot the time evolution of the distance of the average strategy from its equilibrium, $\left\| \textup{avg}(\bs{x}(t)) - \bar{\sigma} \right\|$ (Fig. \ref{fig:avg}) and that of the integral control signal $\sigma(t)$ from its equilibrium $\bar{\sigma}$ (Fig. \ref{fig:sigma}).

In our numerical experience, the control gain $k$ in \eqref{aggregative.dynamics.vector} can sensibly speed up the convergence to an equilibrium, but only up to a critical value above which oscillations may occur. Indeed, the inequality in \eqref{stab.cond} does not hold for large $k$.

\medskip
\section{Conclusion and outlook} \label{sec:conclusion}

In the context of aggregative games, continuous-time integral dynamics with semi-decentralized computation and information exchange can ensure global asymptotic convergence to aggregative equilibria, provided that the control gain is appropriately chosen. 

Future research will focus on extending the convergence analysis to the case of general network games, possibly with coupling constraints, and with relaxed technical assumptions.

\medskip
\section{Acknowledgments}
The authors would like to thank Mr.~G.~Belgioioso (TU Eindhoven),  Mr.~T.W.~Stegink and Mr.~T.~Van Damme (U.~Groningen) for valuable discussions on related topics.

\smallskip
\section*{Appendix: Proofs}

\subsection*{Proof of Lemma \ref{lem:cnes.agg.eq}}
By convexity and the minimum principle, a collective vector $\bs{\bar{x}} = {\rm col}\left(\bar{x}^1, \ldots, \bar{x}^N\right)$ is an aggregative equilibrium exists if and only if 
\[
(z^i-y^i)^\top \left.\frac{\partial J^i}{\partial y^i}(y_i,{\rm avg}(\bs{\bar{x}}))\right|_{y^i=\bar{x}^i} \ge 0, \ \ \forall i\in \mathcal{I}, \ \forall z^i \in \mathcal{X}^i
\]
In view of \eqref{cost.aggregative}, the variational inequality above can be equivalently rewritten as in \eqref{VI}.
\hfill $\blacksquare$

\subsection*{Proof of Lemma \ref{lem:agg.eq.and.dyn.eq}}
Let $\left( \bs{\bar x}, \bar \sigma \right)$ be an equilibrium for \eqref{aggregative.dynamics}, i.e.\, 
\[
\ba{rcll}
0 &=& \Pi_{\mathcal{X}^i} \left(\bar x^i, - \partial f^i\left( \bar x^i \right) -  C \bar\sigma\right), & \forall i\in \mathcal{I}, \\
\smallskip
0 &=& -\bar \sigma + \frac{1}{N} \sum_{j\in \mathcal{I}} \bar x^j
\ea
\]
or, equivalently, 
\[
\ba{rcll}
0 &=& \Pi_{\mathcal{X}^i} \left(\bar x^i \, , \, - \partial f^i\left( \bar x^i \right) -  C {\rm avg}(\bs{\bar x})\right), & \forall i\in \mathcal{I}, \\
\smallskip
\bar \sigma &=& {\rm avg}(\bs{\bar x}).
\ea
\]

If $\bs{\bar x} \in {\rm int}(\bs{\mathcal{X}})$, then  $\Pi_{\mathcal{X}^i} \left( \bar x^i \, , \, - \partial f^i( \bar x^i)-  C {\rm avg}(\bs{\bar x}) \right) = - \partial f^i(\bar x^i)-  C {\rm avg}(\bs{\bar x})$, for all $i\in \mathcal{I}$, hence the equilibrium is such that 
$0 = - {\partial f^i}(\bar x^i) - C {\rm avg}(\bs{\bar x}), \ \forall i\in \mathcal{I}$.
This equation implies that $\bs{ \bar x }$ satisfies the variational inequality in \eqref{VI}, hence, by Lemma \ref{lem:cnes.agg.eq}, it is an aggregative equilibrium. 
 
If $\bs{\bar x} \in {\rm bdry}(\bs{\mathcal{X}})$, 
then there exists a nonempty subset of indices $\mathcal{I}_{{\rm bdry}} \subseteq \mathcal{I}$ such that  $\bar x^i \in {\rm bdry}(\mathcal{X}^i)$ for all $i\in \mathcal{I}_{{\rm bdry}}$. For each $i\in \mathcal{I}_{{\rm bdry}}$, we necessarily have  (see \cite{gadjov:pavel:17cdc}) that
\[
-\partial f^i(\bar x^i) - C{\rm avg}(\bar x) \in N_{\mathcal{X}^i}(\bar x^i)
\]
and by definition of the normal cone operator $N_{\mathcal{X}^i}$, it follows that, for all $i \in \mathcal{I}_{{\rm bdry}}$,
\[
\left( z-\bar x^i\right)^\top \left({\partial f^i}(\bar x^i) + C {\rm avg}(\bs{\bar x}) \right) \ge 0, \ \ \ \forall z\in \mathcal{X}^i.
\]

On the other hand, for each $i\in \mathcal{I}\setminus \mathcal{I}_{{\rm bdry}}$, we have that  
$\Pi_{\mathcal{X}^i} \left( \bar x^i, - \partial f^i(\bar x^i) -  C {\rm avg}(\bs{\bar x}) \right) = - {\partial f^i}(\bar x^i) - C {\rm avg}(\bs{\bar x})$, hence $- \frac{\partial f^i}{\partial x^i}(\bar x^i)-  C {\rm avg}(\bs{\bar x})=0$. We conclude that even when $\bs{\bar x} \in {\rm bdry}(\bs{\mathcal{X}})$, the variational inequality in \eqref{VI} holds and, by Lemma \ref{lem:cnes.agg.eq}, $\bs{\bar x}$ is an aggregative equilibrium.

Conversely, if $\bs{\bar x}$ is an aggregative equilibrium, then $-\partial f^i(\bar x^i) - C{\rm avg}(\bs{\bar x}) \in N_{\mathcal{X}^i}(\bar x^i)$ for all $i\in \mathcal{I}$. Since $\Pi_{\mathcal{X}^i}(\bar x^i, v) = \proj_{T_{\mathcal{X}^i(\bar x^i)}}(v)$ for all $v\in \mathbb{R}^n$, we have that $\Pi_{\mathcal{X}^i}\left(\bar x^i \, , \, -\partial f^i(\bar x^i) - C{\rm avg}(\bar x) \right) = 0$. The thesis then follows by setting $\bar \sigma = {\rm avg}(\bs{\bar x})$.
\hfill $\blacksquare$

\subsection*{Proof of Lemma \ref{cor.key}}
We adapt the approach in \cite{gadjov:pavel:17cdc} to the considered storage function $V$.
Since ${\partial V}(\bs{x}) = \bs{x} - \proj_{\bar{\bs{\mathcal{X}}}}(\bs{x})$, we have that
\begin{multline}
\label{uno}
\frac{\partial V}{\partial x^i}(\bs{x})^\top \Pi_{\mathcal{X}^i} \left( x^i \, , \, - \partial f^i(x^i) + u^i\right)
= \\
\left(x^i - \left[ \proj_{\bar{\bs{\mathcal{X}}}}(\bs{x}) \right]_i \right)^\top \Pi_{\mathcal{X}^i} \left( x^i \, , \, - {\partial f^i}(x^i) + u^i\right)=
\\
-\left(x^i - \left[ \proj_{\bar{\bs{\mathcal{X}}}}(\bs{x}) \right]_i \right)^\top \left[
{\partial f^i}(x^i) - u^i \right. \\ - \Pi_{\mathcal{X}^i} \left( x^i \, , \, - {\partial f^i}(x^i) + u^i \right)
\left. -{\partial f^i}(x^i) + u^i
\right],
\end{multline}
where $\left[ \proj_{\bar{\bs{\mathcal{X}}}}(\bs{x}) \right]_i$ denotes the $i$th block-component of the vector $\proj_{\bar{\bs{\mathcal{X}}}}(\bs{x})$.
Now, as a consequence of Moreau's decomposition theorem, it holds that
\begin{multline}
- \Pi_{\mathcal{X}^i} \left( x^i \, , \, - {\partial f^i}(x^i) + u^i \right)
 -{\partial f^i}(x^i) + u^i \\ 
= \proj_{N_{\mathcal{X}^i}(x^i)}
\left(
- \partial f^i(x^i) + u^i.
\right)
\end{multline}
and, by definition of $N_{\mathcal{X}^i}$ and noting that $\left[ \proj_{\bar{\bs{\mathcal{X}}}}(\bs{x})  \right]_i \in \mathcal{X}^i$, we have that
\begin{multline}
\left( \left[ \proj_{\bar{\bs{\mathcal{X}}}}(\bs{x}) \right]_i - x^i\right)^\top  \\
\left[
- \Pi_{\mathcal{X}^i} \left( x^i \, , \, - {\partial f^i}(x^i) + u^i \right)
- \partial f^i(x^i) + u^i
\right]\le 0.
\end{multline}

From the latter and \eqref{uno}, it then follows that 
\begin{multline}
\label{uno-due}
\frac{\partial V}{\partial x^i}(\bs{x})^\top 
\Pi_{\mathcal{X}^i} \left( x^i \, , \, - {\partial f^i}(x^i) + u^i\right)\\
\le 
-\left( x^i - \left[\proj_{\bar{\bs{\mathcal{X}}}}(\bs{x})\right]_i \right)^\top 
\left( {\partial f^i}(x^i) - u^i
\right).
\end{multline}

Now, we note that since $\bar{\bs{\mathcal{X}}}$ is the set of aggregative equilibria, $\proj_{\bar{\bs{\mathcal{X}}}}(\bs{x})$ is an aggregative equilibrium, and by Lemma \ref{lem:agg.eq.and.dyn.eq},  the pair  $\left( \proj_{\bar{\bs{\mathcal{X}}}}(\bs{x}), {\rm avg}\left( \proj_{\bar{\bs{\mathcal{X}}}}(\bs{x}) \right) \right)$ is an equilibrium for the dynamics in  \eqref{aggregative.dynamics}, hence it satisfies 
\be\label{eq.star}
\bs{0} = \Pi_{\bar{\bs{\mathcal{X}}}}\left( \proj_{\bar{\bs{\mathcal{X}}}}(\bs{x}) \, , \, - {\partial f}\left( \proj_{\bar{\bs{\mathcal{X}}}}(\bs{x}) \right) + \bar{\bs{u}} \right)
\ee
where $\bar{\bs{u}} = -\bs{1}_N \otimes \left( C \, {\rm avg}\left(\proj_{\bar{\bs{\mathcal{X}}}}(\bs{x}) \right) \right)$.
Again by Moreau's decomposition theorem, the latter yields
\begin{multline}
- {\partial f^i}\left( \left[\proj_{\bar{\bs{\mathcal{X}}}}(\bs{x})\right]_i \right) + \bar{u}^i \\ 
=  \proj_{N_{\mathcal{X}^i}(\left[\proj_{\bar{\bs{\mathcal{X}}}}(\bs{x})\right]_i)}
\left(-{\partial f^i}( \left[\proj_{\bar{\bs{\mathcal{X}}}}(\bs{x})\right]_i ) + \bar{u}^i \right),
\end{multline}
which, by definition of $N_{\mathcal{X}^i}\left(\left[\proj_{\bar{\bs{\mathcal{X}}}}(\bs{x})\right]_i\right)$,  implies that
\be\label{due}
\left( -{\partial f^i}([ \proj_{ \bar{\bs{\mathcal{X}}} }(\bs{x}) ]_i)+ u_i^*\right)^\top \left( x^i - \left[\proj_{\bar{\bs{\mathcal{X}}}}(\bs{x})\right]_i \right) \le 0.
\ee
We conclude that 
\begin{multline}
\frac{\partial V}{\partial x^i}(\bs{x})^\top 
\Pi_{\mathcal{X}^i}\left( x^i \, , \, -\partial f^i(x^i)+ u^i\right) = \\
\le -\left( x^i - \left[\proj_{\bar{\bs{\mathcal{X}}}}(\bs{x})\right]_i \right)^\top 
\left(
\partial f^i(x^i) - u^i
\right) \\
\le \left( x^i - \left[\proj_{\bar{\bs{\mathcal{X}}}}(\bs{x})\right]_i \right)^\top 
\left( -{\partial f^i}(x^i) + u^i \right)  \\
+
\left({\partial f^i}\left(\left[\proj_{\bar{\bs{\mathcal{X}}}}(\bs{x})\right]_i\right) - \bar{u}^i  \right)^\top \left( x^i - \left[\proj_{\bar{\bs{\mathcal{X}}}}(\bs{x})\right]_i \right),
\end{multline}
where the first equality holds by \eqref{uno-due}, and the second one holds in view of \eqref{due}. 
The thesis then follows by summing up the inequalities over $i\in \mathcal{I}$.
\hfill $\blacksquare$

\subsection*{Proof of Proposition \ref{prop:lyapunov.decrease}}

For ease of notation, let us define $\bar{\sigma} := {\rm avg}\left( \proj_{\bar{\bs{\mathcal{X}}}}(\bs{x}) \right)$. We first note that
\begin{multline} \label{Lyapunov.in.constr}
\dot W(\bs{x}, \sigma) = \\ 
\begin{bmatrix} \frac{\partial W}{\partial \bs{x}}( \bs{x}, \sigma ) 
\smallskip 
\\  
\frac{\partial W}{\partial \sigma}(\bs{x}, \sigma) \end{bmatrix}^\top
\begin{bmatrix}
\Pi_{\bs{\mathcal{X}}} \left( \bs{x} \, , \, -\partial f(\bs{x})-\bs{1}_N \otimes (C \sigma) \right)\\
k \left( -\sigma + {\rm avg}(\bs{x}) \right)
\end{bmatrix} \\
=
-\left(\bs{x} - \proj_{\bar{\bs{\mathcal{X}}}}(\bs{x}) \right)^\top \left( {\partial f}(\bs{x})-{\partial f}(\proj_{\bar{\bs{\mathcal{X}}}}(\bs{x})) \right) \\
-(x-\proj_{\bar{\bs{\mathcal{X}}}}(\bs{x}))^\top (\bs{1}_N\otimes ( C(\sigma - \bar{\sigma} ) ) ) \\
+\left(\sigma-\bar{\sigma}\right)^\top k (-\sigma + \bar{\sigma} + {\rm avg}(\bs{x})-{\rm avg}( \proj_{\bar{\bs{\mathcal{X}}}}(\bs{x}) )), 
\end{multline}
where we have used Lemma \ref{cor.key} with $\bs{u} = -\bs{1}_N \otimes (C\sigma)$.

Since the functions $\{f^i\}_{i \in \mathcal{I}}$ are twice continuous differentiable and the sets $\mathcal{X}^i$'s are convex, we have that 
\begin{multline}
-\left(\bs{x} - \proj_{\bar{\bs{\mathcal{X}}}}(\bs{x})\right)^\top \left( {\partial f}(\bs{x})-{\partial f}\left( \proj_{\bar{\bs{\mathcal{X}}}}(\bs{x}) \right) \right) \\ 
=
-\textstyle \sum_{i\in \mathcal{I}} (x^i - \left[ \proj_{\bar{\bs{\mathcal{X}}}}(\bs{x}) \right]_i)^\top \qquad \qquad \\
\left[ {\partial^2 f^i}(\hat x^i) \right] \left( x^i - \left[ \proj_{\bar{\bs{\mathcal{X}}}}(\bs{x}) \right]_i \right),
\end{multline}
for some $\hat x^i \in \mathcal{X}^i$, $i\in \mathcal{I}$. 
By $\ell$-strong convexity of $f^i$ on $\mathcal{X}^i$, i.e., $\partial^2 f^i(\xi) \succcurlyeq \ell I_n$ for all $\xi \in \mathcal{X}^i$, we derive that 
\begin{multline}
\label{quadratic.term.x}
-\left( \bs{x} - \proj_{\bar{\bs{\mathcal{X}}}}(\bs{x}) \right)^\top \left( {\partial f}(\bs{x})-{\partial f}( \proj_{\bar{\bs{\mathcal{X}}}}(\bs{x}) ) \right) \\ 
\le -\dst\sum_{i\in \mathcal{I}} \ell \left\| x^i - \left[ \proj_{\bar{\bs{\mathcal{X}}}}(\bs{x}) \right]_i \right\|^2.
\end{multline}
Furthermore, we have that 
\begin{multline}
\label{cross-term}
\left( \bs{x} - \proj_{\bar{\bs{\mathcal{X}}}}(\bs{x}) \right)^\top \left( \bs{1}_N\otimes ( C (\sigma- \bar{\sigma}  )) \right)
\\ = \dst\sum_{i\in \mathcal{I}} \left( x^i - \left[ \proj_{\bar{\bs{\mathcal{X}}}}(\bs{x}) \right]_i \right)^\top C \left( \sigma-\bar \sigma \right)
\end{multline}
and
\begin{multline}
\label{cross-term-2}
\left(\sigma - \bar\sigma \right)^\top k \left( {\rm avg}(\bs{x}) - {\rm avg}( \proj_{\bar{\bs{\mathcal{X}}}}(\bs{x}) ) \right)
\\ = \left( \sigma-\bar\sigma \right)^\top k  \dst \tfrac{1}{N} \sum_{i\in \mathcal{I}} x^i -  \left[ \proj_{\bar{\bs{\mathcal{X}}}}(\bs{x}) \right]_i.
\end{multline}

By rearranging the terms in matrix form, we obtain that 
$$\dot{W}( \bs{x}, \sigma ) \leq -\left\| 
\left[
\begin{matrix}
\bs{x} - \proj_{\bar{\bs{\mathcal{X}}}}(\bs{x}) \\
\sigma - \textrm{avg}( \proj_{\bar{\bs{\mathcal{X}}}}(\bs{x}) )
\end{matrix}
\right]
\right\|_M^2\,,$$
with matrix 
\[
M := \begin{bmatrix}
\ell \, I_{nN} & -\frac{1}{2}  \begin{bmatrix} C+\frac{k}{N} I_n \\ \vdots \\ C+\frac{k}{N} I_n\end{bmatrix} 
\smallskip\\
\star &  k \, I_{nN}
\end{bmatrix}.
\]
Finally, we apply the Gershgorin theorem to the matrix $M$. Namely, the condition in \eqref{stab.cond} implies that 
\[
\min \{ \ell, k\} > \tfrac{1}{2} \max_{i\in \mathcal{I}} \sum_{j\ne i}^N |c_{i,j}|+ |c_{i,i}+\tfrac{k}{N}|, 
\]
which in turn implies that $M \succ 0$. The proof then follows with $\epsilon := \lambda_{\min}(M) > 0$.
\hfill $\blacksquare$


\subsection*{Proof of Theorem \ref{main.aggreg}}
The dynamics in \eqref{aggregative.dynamics.vector} represent a projected dynamical system with discontinuous right-hand side \cite{nagurney}. Thus, its solutions must be intended in a Carath\`{e}odory sense. 

By \cite[Th. 2.5]{nagurney} and \cite[Prop. 2.2]{cherukuri.scl16}, since the mapping ${\rm col}\left( \bs{x}, \sigma \right) \mapsto {\rm col}\left( - \partial f( \bs{x} ) - \bs{1}_N \otimes ( C \sigma ) \, , \, -\sigma +{\rm avg}( \bs{x} ) \right)$ is Lipschitz continuous on the closed convex set $\boldsymbol{\mathcal{X}}$, for any initial condition $(\boldsymbol{x}_0, \sigma_0) \in \boldsymbol{\mathcal{X}} \times \R^n$,  
there exists a unique Carath\`{e}odory solution to \eqref{aggregative.dynamics.vector}, $\left(\boldsymbol{x}(\cdot), \sigma(\cdot) \right)$,  
such that $\left(\boldsymbol{x}(0), \sigma(0) \right) = \left(\boldsymbol{x}_0, \sigma_0 \right)$, $\left( \boldsymbol{x}(t), \sigma(t) \right) \in \boldsymbol{\mathcal{X}} \times \R^n$ for all $t\in [0, +\infty)$, and that is uniformly continuous with respect to the initial condition. 
Let $\mathcal{S}_0$ be a compact sublevel set of the Lyapunov function $W$ that contains the initial condition ${\rm col}\left( \bs{x}_0, \sigma_0 \right)$. Then, the intersection set $\mathcal{C} := \mathcal{S}_0 \cap \left( \bs{\mathcal{X}} \times \mathbb{R}^n \right)$ is a compact set. 
Let us consider the Carath\`{e}odory solution $(\bs{x}(t), \sigma(t))$ issuing from ${\rm col}\left( \bs{x}_0, \sigma_0 \right)$ and evaluate the Lyapunov function $W(\bs{x},\sigma)$ along such solution. By \eqref{strict.lyap}, and by a chain rule for absolutely continuous functions, 
for almost every $t\ge 0$,  
\[
\dot W (\bs{x}(t),\sigma(t)) \le - \lambda_{\min}(M) \, W ( \bs{x}(t) , \sigma(t))\le 0.
\]  
Thus, $W$ is non-increasing along $(\bs{x}(t), \sigma(t))$, 
and we conclude that $\mathcal{C}$ is a forward invariant set. Using a Gronwall--Bellman inequality,  see the final arguments in the proof of \cite[Th. 2]{brogliato:goeleven:05}, we obtain that 
\[
W (\bs{x}(t),\sigma(t)) \le {\rm exp}\left(-\lambda_{\min}(M) t \right) \, W (\bs{x}_0, \sigma_0 ),
\]  
and therefore 
\begin{multline}
\left\| \begin{bmatrix}
\bs{x}(t) - \proj_{\bar{\bs{\mathcal{X}}}}(\bs{x}(t))\\
\sigma(t) - {\rm avg}\left( \proj_{\bar{\bs{\mathcal{X}}}}(\bs{x}(t)) \right)
\end{bmatrix} \right\|
\\
\le {\rm exp}\left(- \tfrac{\lambda_{\min}(M)}{2} t \right) 
\left\| 
\begin{bmatrix}
\bs{x}_0 - \proj_{\bar{\bs{\mathcal{X}}}}(\bs{x}_0)\\
\sigma_0 - {\rm avg}\left( \proj_{\bar{\bs{\mathcal{X}}}}(\bs{x}_0) \right)
\end{bmatrix} \right\|,
\end{multline}
which proves global exponential convergence.
\hfill $\blacksquare$

\balance

\medskip
\bibliographystyle{IEEEtran}
\bibliography{library,int-nask-seeking}

\end{document}